\documentclass{amsart}
\usepackage{calc}
\usepackage{color}

\usepackage{enumerate}

\usepackage{amsmath}
\usepackage{amsfonts}
\usepackage{amssymb}
\usepackage{amsxtra}
\usepackage{amsthm}
\usepackage{array}

\usepackage{graphics}
\usepackage{graphicx}
\usepackage{epic}
\usepackage{eepic}

\usepackage{supertabular}

\def\doprivate{01}

\newcommand{\changed}[1]{#1}

\newcommand{\R}{\mathbb{R}}
\newcommand{\RC}{\R^{d+1}_+}
\newcommand{\Rd}{\R^d}
\newcommand{\defm}[1]{{\it #1}}
\newcommand{\loc}{\operatorname{loc}}
\newcommand{\supp}{\operatorname{supp}}

\theoremstyle{definition}
\newtheorem{definition}{Definition}
\theoremstyle{plain}
\newtheorem{theorem}{Theorem}

\newtheorem{lemma}{Lemma}
\newtheorem{proposition}{Proposition}

\theoremstyle{remark}

\begin{document}

%\thanks{This material is based upon work supported by an SAP/Stanford Graduate Fellowship and by the National Science Foundation under Grant no. DMS 0104019.}
\thanks{This material is based upon work supported by an SAP/Stanford Graduate Fellowship and by the National Science Foundation under
Grant no. DMS 0104019. Any opinions, findings, and conclusions or recommendations expressed in this material are those of the author and
do not necessarily reflect the views of the National Science Foundation.}


%%\begin{frontmatter}

%\title[Num. investigation of well-posedness of systems of c. laws]{Numerical investigation of the well-posedness of systems of conservation laws}
%%\title{A possible counterexample to uniqueness of entropy solutions and Godunov scheme convergence\thanksref{blah}}
\title[Wellposedness of entropy solutions]{A possible counterexample to wellposedness of entropy solutions and to Godunov scheme convergence}
%Numerical investigation of the well-posedness of systems of conservation laws}
%%\author{Volker Elling\corauthref{blub}}
%%\thanks[blah]{This material is based upon work supported by an SAP/Stanford Graduate Fellowship and by the National Science Foundation under Grant no. DMS 0104019. Any opinions, findings, and conclusions or recommendations expressed in this material are those of the author and do not necessarily reflect the views of the National Science Foundation.}
\author{Volker Elling}
\address{Brown University\\Division of Applied Mathematics\\182 George Street\\Providence, RI 02912}
\email{velling@stanfordalumni.org}
\urladdr{http://www.dam.brown.edu/people/volker}
%\author{Volker Elling\thanks{SCCM Program, Gates Building 2B, Stanford University, Stanford, CA 94305-9025 ({\tt velling@stanford.edu}). This material is based upon work supported by an SAP/Stanford Graduate Fellowship and by the National Science Foundation under Grant no. DMS 0104019. Any opinions, findings, and conclusions or recommendations expressed in this material are those of the author and do not necessarily reflect the views of the National Science Foundation.}}
%\date{March 22, 2003}

%%\address[blub]{SCCM Program\\Gates Building 2B\\Stanford University\\Stanford, CA 94305-9025}
%%\ead{velling@stanford.edu}
%%\ead[url]{http://www.stanford.edu/\~{}velling}
%\email{velling@stanford.edu}
%\address{Volker Elling\\SCCM Program\\Gates Building 2B\\Stanford University\\Stanford, CA 94305-9025}
%\urladdr{http://www.stanford.edu/\~{}velling}

\begin{abstract}
A particular case of initial data for the two-dimensional Euler equations is studied numerically. 
The results show that the Godunov method does not always converge to the physical solution, at least not on 
feasible grids. 
Moreover, they suggest that entropy solutions (in the weak entropy inequality sense) are not well-posed. 
\end{abstract}

\maketitle
%%\end{frontmatter}

%%\begin{keywords}
%%Conservation law, wellposedness, entropy solution, Riemann problem, shock, contact discontinuity, compressible Euler equations,
%%entropy-entropy flux pair
%%\end{keywords}

%%\begin{AMS}
%%35L65, 35L67, 76L05, 76H05, 76N10
%%\end{AMS}
%%pick one primary; also check: 35M, 76H05, 76J20, 76K05, 76M12, 76N15, 76N10

%\parindent=0cm%
%\parskip=.3cm%

\section{Introduction}

%%\subsection{Systems of conservation laws}

%%\UCON{cut down the introductory stuff a bit. not all needed for purpose of paper. or make more coherent.}

%%\UCON{for jcp, explain in more detail how people expected to use discrete entropy inequality for showing that the solution is physical. put REFS in here!}

Consider the Cauchy problem for a system of hyperbolic conservation laws,
\begin{alignat}{1}
    \frac{\partial u}{\partial t} + \nabla\cdot(\Vec f(u)) &= 0, \label{eq:system} \\
    u(0,\cdot) &= u_0, \label{eq:system-inicond}
\end{alignat}
where $u=u(t,\Vec x):\RC:=(0,\infty)\times\R^d\rightarrow P\subset\R^m$ is the desired solution ($P$ the set of physically reasonable values),
$\Vec f=(f^i)$, $f^i:P\rightarrow\R^m$, the (smooth) \defm{flux} function,
$u_0:\R^d\rightarrow P$ initial data. Here and in the sequel ``$\nabla,\Delta,\cdot$'' are meant with respect to $\Vec x$.

%%For $d=1$, the system is called \defm{strictly hyperbolic} if
%%$\frac{\partial f}{\partial u}(w)$ has $m$ distinct real eigenvalues (for all $w\in P$); it is called \defm{strongly hyperbolic}
%%if the eigenvalues are real and semisimple (i.e.\ if
%%there is a complete set of eigenvectors). One says that the $i$th field is \defm{linearly degenerate}
%%if $\frac{\partial\lambda^i}{\partial u}(w)r^i(w)=0$ for all $w\in P$, \defm{genuinely nonlinear} if $\frac{\partial\lambda^i}{\partial u}(w)r^i(w)\neq 0$
%%for all $w\in P$ (here, $r^i(w)$ is the right eigenvector of $\frac{\partial f}{\partial u}(w)$ corresponding to $\lambda^i(w)$).
%%For $d>1$, these terms are used if the conditions hold for the 1D system with flux $\Vec f\cdot\Vec n$, for all unit vectors $\Vec n$.

An important example of hyperbolic systems of conservation laws are the (nonisentropic) compressible Euler equations:
\begin{alignat}{1}
    \rho_t + \nabla\cdot(\rho\Vec v) &= 0, \notag \\
    (\rho v^i)_t + \nabla\cdot(\rho v^i\Vec v) + p_{x_i} &= 0 \qquad (i=1,\dotsc,d), \notag \\
    (\rho e)_t + \nabla\cdot((\rho e+p)\Vec v) &= 0. \label{eq:nonisentropic-euler}
\end{alignat}
Here, $\rho$ is density, $\Vec v=(v^i)$ velocity, $e$ specific energy, which decomposes into
\begin{alignat}{1}
    e &= \frac{|\Vec v|^2}{2}+q; \label{eq:specific-energy}
\end{alignat}
the first summand is specific kinetic energy, $q$ is specific internal energy.
The pressure is a function of $\rho,q$; a common choice is the polytropic pressure law
\begin{alignat}{1}
    p &= (\gamma-1)\rho q \label{eq:nonisentropic-polytropic-pressure}
\end{alignat}
($1<\gamma\leq\frac{5}{3}$; for air, $\gamma=\frac{7}{5}$).
The set of admissible values is
\[ P=\{q>0,\ \rho>0\}. \]

It is well-known that \eqref{eq:system} and \eqref{eq:system-inicond} need not have a global smooth solution, even if the initial data $u_0$ is smooth.
For this reason, one has to study \defm{weak solutions}, defined as functions $u\in L^1_{\loc}(\RC;P)$ that satisfy
\begin{alignat}{1}
    -\int_0^\infty\int_{\R^d} u\frac{\partial\phi}{\partial t}+\Vec f(u)\cdot\nabla\phi~d\Vec x~dt &= \int_{\R^d}u_0(\Vec x)\phi(0,\Vec x)~d\Vec x, \label{eq:weak}
\end{alignat}
for all test functions $\phi\in C^\infty_c(\overline{\RC})$. Moreover, there can be more than one weak solution, so it is
necessary to impose an additional condition, called \defm{entropy condition}, to single out a unique weak solution (the \defm{entropy solution}).

One definition of entropy solutions is the \defm{vanishing viscosity} (VV) definition; it requires that $u$ is the limit of the sequence
$(u^\epsilon)_{\epsilon>0}$ of solutions of
\begin{alignat}{1}
    \frac{\partial u^\epsilon}{\partial t} + \nabla\cdot(\Vec f(u^\epsilon)) &= \epsilon\Delta u^\epsilon\qquad\text{in $\RC$,}\label{eq:viscous-system} \\
    u^\epsilon(0,\cdot) &= u_0\qquad\text{on $\{0\}\times\R^d$}. \label{eq:viscous-system-inicond}
\end{alignat}
The limit is taken in some suitable topology, usually as a boundedly almost everywhere limit. We call such a function $u$ a \defm{VV solution}.

Another definition uses \defm{entropy/entropy flux} (EEF) pairs $(\eta,\Vec\psi)$,
where $\eta:P\rightarrow\R$ is a smooth strictly convex function, called \defm{entropy}, whereas
$\Vec\psi=(\psi^1,\dotsc,\psi^d)'$ with smooth $\psi^i:P\rightarrow\R$ is called \defm{entropy flux}; $\eta$ and $\Vec\psi$ are required to satisfy
\begin{alignat}{1}
    \frac{\partial\psi^i}{\partial u^\alpha}&=\sum_{\beta=1}^m\frac{\partial\eta}{\partial u^\beta}\frac{\partial f^{i\beta}}{\partial u^\alpha}
    \qquad(i=1,\dotsc,d,\ \alpha=1,\dotsc,m).
    \label{eq:eef}
\end{alignat}
By multiplying \eqref{eq:viscous-system} from the left with $\eta'(u^\epsilon)$ and using \eqref{eq:eef},
one obtains
\begin{alignat}{1}
    \frac{\partial(\eta\circ u^\epsilon)}{\partial t} + \sum_{i=1}^d\frac{\partial(\psi^i\circ u^\epsilon)}{\partial x^i} &= \epsilon\Delta(\eta\circ u^\epsilon) - \epsilon\sum_{i=1}^d\eta''(u^\epsilon)\frac{\partial u^\epsilon}{\partial x^i}\frac{u^\epsilon}{\partial x^i} 
    \leq \epsilon\Delta(\eta\circ u^\epsilon). \label{eq:weak-viscous}
\end{alignat}
(here, we used that $\eta$ is convex).
Upon multiplying the last equation with a \emph{nonnegative} test function $\phi$ and integrating by parts, this yields
\begin{alignat}{1}
    -\int_0^\infty\int_{\R^d}\eta(u^\epsilon)\frac{\partial\phi}{\partial t} + \Vec\psi(u^\epsilon)\cdot\nabla\phi~d\Vec x~dt &\leq \epsilon\int_0^\infty\int_{\R^d}\eta(u^\epsilon)\Delta\phi~d\Vec x~dt + \int_{\Rd}\eta(u_0)~d\Vec x. \label{eq:wwv}
\end{alignat}
If, as assumed above, $(u^\epsilon)\rightarrow u$ boundedly almost everywhere, then \eqref{eq:wwv} implies
\begin{alignat}{1}
    -\int_0^\infty\int_{\R^d}\eta(u)\frac{\partial\phi}{\partial t} + \Vec\psi(u)\cdot\nabla\phi~d\Vec x~dt & \leq \int_{\R^d}\eta(u_0)~d\Vec x. \label{eq:eef-entropy}
\end{alignat}
Functions $u$ that satisfy \eqref{eq:eef-entropy} for \emph{all} EEF flux pairs are called \defm{EEF solutions} (of \eqref{eq:system}). 
As we have shown, VV solutions are necessarily EEF solutions.

In the literature, the term \defm{entropy solution} is used to refer either to EEF or to VV solutions, often without explicit mention, 
because it has been assumed that the two definitions are equivalent for the Euler equations and many other 
physically relevant systems (see \cite{serre-1} p.\ 101, \cite{dafermos-book} p.\ 49, \cite{godlewski-raviart} p.\ 32; see the discussion 
in Section \ref{section:related-work} for verified special cases). However, for the purposes of this paper it is necessary
to distinguish the two notions, as we will discuss a possible numerical counterexample to their equivalence.

The (gas-dynamic) specific entropy $s$ is defined as
\begin{alignat}{1}
    s &= \log q+(1-\gamma)\log\rho; \label{eq:gas-dynamic-entropy}
\end{alignat}
\begin{alignat}{1}
    \eta &:= -\rho s,\quad \psi^i:=-\rho sv^i
\end{alignat}
provides an EEF pair for the Euler equations. 

A common simplification is to assume that $s$ is constant in space and time. This yields the \emph{isentropic} Euler equations
\begin{alignat}{1}
     \rho_t + \nabla\cdot(\rho\Vec v) &= 0, \notag \\
    (\rho v^i)_t + \nabla\cdot(\rho v^i\Vec v) + p_{x_i} &= 0 \qquad (i=1,\dotsc,d) \label{eq:isentropic-euler}
\end{alignat}
with
\begin{alignat}{1}
    p(\rho) &= \rho^\gamma. \label{eq:isentropic-polytropic-pressure}
\end{alignat}
In this case, $P=\{\rho>0\}$. An EEF pair is provided by the specific energy $e$,
\begin{alignat*}{1}
    e &= \frac{|\Vec v|^2}{2}+\frac{\rho^{\gamma-1}}{\gamma-1},
\end{alignat*}
with
\begin{alignat*}{1}
    \eta &:= \rho e, \quad \psi^i:=(\rho e+p)v^i.
\end{alignat*}

It is cumbersome to verify the EEF condition \eqref{eq:eef-entropy} directly, not to mention the VV condition. 
There are easier criteria for piecewise smooth functions, which we define in the following customized way:

\begin{definition}
    \begin{enumerate}
    \item
        A point $(t,\Vec x)\in\RC$ is called \defm{point of smoothness} if $u$ is $C^\infty$ in a small neighbourhood of $(t,\Vec x)$.
    \item
        A point $(t,\Vec x)\in\RC$ is called \defm{point of piecewise smoothness} of $u$ if 
        there is a $C^\infty$ diffeomorphism $\Phi$ of a ball $V$ around $0$ in $\R^{d+1}$ onto a neighbourhood
        $B$ of $(t,\Vec x)=\Phi(0)$  so that $u\circ\Phi$ is $C^\infty$ on $B_-$ and on $B_+$ (where $B_\pm:=\Phi(V_\pm)$, $V_\pm:=\{y\in V:y_1\gtrless 0\})$;
        for later use, let $S$ be the surface $\Phi(V\cap(\{0\}\times\R^d))$, $n=(n^t,\Vec n)\in\R^{d+1}$
        a unit normal to $S$ in $(t,\Vec x)$ pointing into $B_+$;
        let $u_+,u_-$ be the one-sided limits of $u$ in $(t,\Vec x)$ within $B_-$ resp.\ $B_+$). We also require $\vec n\neq 0$.
        %%    For $\Vec n\neq 0$, we define the shock speed 
        %%    \[ \Vec\sigma:=\frac{\Vec n}{n_t},\qquad\sigma:=|\Vec\sigma|. \]
        %%    (Note: points of smoothness are points of piecewise smoothness via an arbitrary, unimportant choice of $\Phi$.)
    \item
        $u$ is called \defm{piecewise smooth} if there is a set $N$ of $d$-dimensional Hausdorff measure $0$ so that
        all points in $\RC-N$ are points of piecewise smoothness.
    \end{enumerate}
\end{definition}

\begin{proposition}
    Let $u$ be piecewise smooth.
    $u$ is an EEF solution of \eqref{eq:system} if and only if 
    \begin{enumerate}
    \item it is a (classical) solution of \eqref{eq:system} in each point of smoothness,
    \item $u(t,\cdot)\rightarrow u_0$ in $L^1_{\loc}$ as $t\downarrow 0$, and
    \item in each point $(t,x)$ of piecewise smoothness it satisfies the Rankine-Hugoniot conditions 
        \begin{alignat}{1}
            (u_+-u_-)n^t + (\Vec f(u_+)-\Vec f(u_-))\cdot\Vec n &= 0 \label{RH}
        \end{alignat}
        and (for all EEF pairs $(\eta,\Vec\psi)$)
        \begin{alignat}{1}
            (\eta(u_+)-\eta(u_-))n^t + (\Vec\psi(u_+)-\Vec\psi(u_-))\cdot\Vec n &\leq 0. \label{RH-EEF}
        \end{alignat}
    \end{enumerate}
    \label{th:piecewise}
\end{proposition}

Proposition \ref{th:piecewise} is well-known (see, for example, Section 11.1.1 in \cite{evans}), as is the following property:
\begin{proposition}
    \label{rem:shock-vn}%
    For the Euler equations \eqref{eq:nonisentropic-euler} resp.\ \eqref{eq:isentropic-euler} (with polytropic gas law 
    \eqref{eq:nonisentropic-polytropic-pressure} resp.\ \eqref{eq:isentropic-polytropic-pressure}),
    \eqref{RH-EEF} is equivalent to the simpler condition that the normal velocity does not increase across discontinuities:
    \begin{alignat*}{1}
        (\vec v_+-\vec v_-)\cdot\vec n &\leq 0.
    \end{alignat*}
\end{proposition}

%%\subsection{Compressible Euler equations}

%%\label{section:euler}

%%The Euler equations are strictly hyperbolic in one dimension; in two or more dimensions, they are strongly, but not strictly, hyperbolic. 
%%They have two genuinely nonlinear and $d$ (nonisentropic) resp.\ $d-1$ (isentropic) linearly degenerate fields.
The Cauchy problem for the Euler equations has several important symmetry properties, including the following:
\begin{proposition}
    \label{th:symmetries}%
    Let $u=(\rho,\Vec v',q)'$ be a weak solution for initial data $u_0=(\rho_0,\Vec v_0',\Vec q_0)'$.
    \begin{enumerate}
        %%\item Galilean invariance:
%\begin{enumerate}
        %%\item Translation invariance: for all $\Vec y\in\R^d$, $u(\cdot+\Vec y,\cdot)$ is a weak solution for initial data $u_0(\cdot+\Vec y)$.
        %%\item Rotation invariance/mirror symmetry:
        %%    For all orthogonal matrices $Q$, $(\rho(Q\Vec x),(Q^*v(Q\Vec x))',q(Q\Vec x))'$ is a weak solution for initial data $u_0(Q\cdot)$.
        %%\end{enumerate}
    \item Change of inertial frame:
        For all $\Vec w\in\R^d$, $(\rho(x+\Vec wt),(\Vec v(x+\Vec wt,t)-\Vec w)',q(x+\Vec wt))'$ is a weak solution for the same initial data $u_0$.
    \item Self-similarity: a function $f:\RC\rightarrow\R^m$ is called \defm{self-similar} if $f(rt,r\Vec x)=f(t,\Vec x)$ for all $r>0$;
        same for functions on $\Rd$.
        If the initial data is self-similar, then for any $r>0$, $u(r\Vec x,rt)$ is a weak solution for the same initial data $u_0$.
    \end{enumerate}
    These symmetries remain true after replacing ``weak'' by ``VV'' or ``EEF''. Analogous symmetries hold for the isentropic case.
\end{proposition}

\section{Example and numerical results}

\label{section:numerical-results}

%(TODO: explain somewhere that it is also important that the second solution coincides with the boundary values in a ball $\{(t,x):|x|>rt\}$.)

\begin{figure}
\input{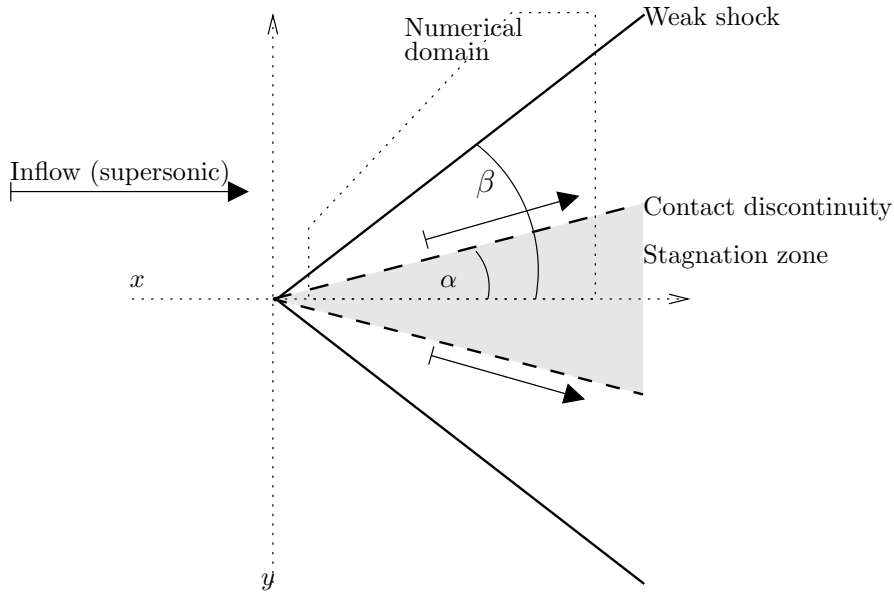}
\caption{Solution T: this initial data is also a steady and self-similar solution 
for the compressible Euler equations in 2D}
\label{fig:example}
\end{figure}

Consider the following set $u_0$ of initial data for \eqref{eq:nonisentropic-euler} with $d=2$ (see Figure \ref{fig:example}):
the data is symmetric under reflection across the $x$-axis and constant in each of four cones centered in the origin (in particular,
constant along rays starting in the origin).
In the origin, two shocks emanate into the first and fourth quadrant;
the area on the left is supersonic inflow (parallel to the $x$-axis); the two areas on the other side of the shocks are denser and hotter
gas, moving parallel to the contact discontinuities (see \cite{courant-friedrichs} Chapter IV C on choosing pre- and 
post-shock values that satisfy the Rankine-Hugoniot conditions; we choose the ones that yield the weaker shock). 
The gas in the stagnation area (enclosed by the contact discontinuities) has the same pressure 
as the post-shock gas on the other side, but velocity $\Vec v=0$.
It is easy to check, using Propositions \ref{th:piecewise} and \ref{rem:shock-vn}, that the steady solution $u(t,\Vec x)=u_0(\Vec x)$ is an 
EEF solution of \eqref{eq:nonisentropic-euler} resp.\ \eqref{eq:isentropic-euler}. Henceforth we refer to it as \defm{Solution T} (for \defm{theoretical}).

\begin{figure}
\includegraphics[angle=90,width=\textwidth]{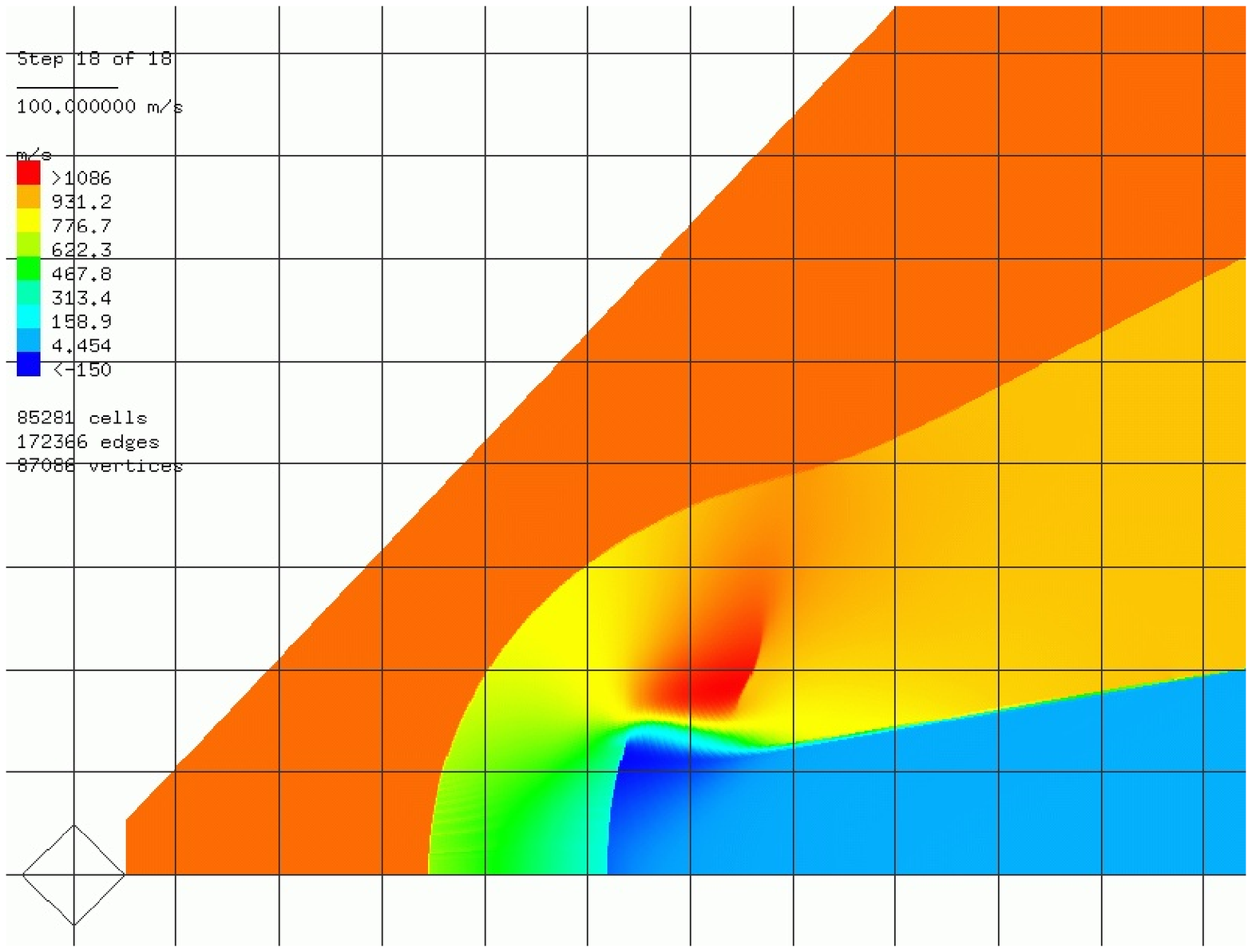}
\caption{%
Solution N:
(Rotate clockwise by $90^o$ to align with dotted area in Figure \ref{fig:example}.)
Each square in the coordinate grid corresponds to a $100~m/s\times100~m/s$ square in the $\vec\xi$ plane. 
Origin marked by diamond (lower left corner). 
Plotted: horizontal velocity.
Godunov scheme for isentropic Euler equations; 
data: $\gamma=1.4$, $\alpha=10^o$; inflow: $\rho=1.19~kg/m^3$, $v=1000~m/s$, $T=20^oC$.
The solution differs significantly from Figure \ref{fig:example}. 
Results for nonisentropic Euler equations
or other numerical schemes are similar.
}
\label{fig:second-vx}
\end{figure}
\changed{However, instead of Solution T, numerical calculations produce the markedly different result in Figure \ref{fig:second-vx} which we call \defm{Solution N} (for \defm{numerical}) in the sequel (of course it is not known to be an exact solution). 
The numerical domain in Figure \ref{fig:second-vx} is indicated as the dotted quadrilateral in Figure \ref{fig:example}. Figure \ref{fig:second-vx} was computed as follows:}
%%\begin{figure}
%%    \epsfig{file=selfsim-rho.epsi,height=.85\textheight}
%%%    \epsfig{file=selfsim-rho-bw.epsi,height=.85\textheight}
%%    \caption{%
%%        Plot of $\rho$, for the nonisentropic Euler equations, computed with the ENO-RF scheme. In contrast to the isentropic
%%        case, there is an extra entropy-jump type contact discontinuity between the two shocks that meet the $\xi_1$ axis.
%%    }
%%\label{fig:second-rho}
%%\end{figure}
adaptive refinement was used to achieve better resolution at same computational cost. To reduce numerical viscosity
the grid was chosen so that near the right domain boundary the edges are aligned with the contact discontinuity and the shock.
In order to capture self-similarity, the computations were done for a grid with moving vertices with coordinates $\Vec x=t\Vec\xi$
($\Vec\xi$ has the dimension of a velocity; its components are called \defm{similarity coordinates}).
The moving-edge modifications discussed in \cite[Section 2.1.6]{elling-diplom} and \cite[Chapter 4]{phd-thesis} were used (the essential idea is to compute numerical fluxes 
across a moving edge by transforming to a steady edge, using invariance under change of inertial frame (Proposition \ref{th:symmetries}),
and to apply an arbitrary approximate Riemann solver to the transformed problem). 
The domain boundaries were chosen so that small perturbations on them propagate into the domain ($\Vec\xi\cdot\Vec n$ ($\Vec n$ outer unit normal) 
in each boundary point is larger than the maximum of $|\Vec v|+c$ in the domain). This allows to prescribe all components of the 
fluxes on the boundary.

Experiments with various modifications were made: changing the numerical scheme (the experiments were repeated for the Godunov scheme \cite{godunov}, the 
Osher-Solomon scheme \cite{solomon-osher}, the ENO-RF scheme \cite{shu-osher-llf}, and a second-order MUSCL code based on the first-order ENO-RF scheme), 
adding more numerical dissipation, refining uniformly rather than adaptively, using a Cartesian grid including origin and lower half-plane, 
or calculating in space rather than similarity coordinates. None of these modifications change the numerical results significantly; in all cases,
the numerical results converge to the same Solution N.

\changed{Solution N appears to be self-similar (i.e.\ steady in similarity coordinates), but it is strongly \emph{unsteady}, so it
is clearly different from Solution T.}

\section{Conclusions about numerical methods}

\label{section:numfail}

\changed{While the discrepancy between Solutions N and T opens many new problems, we can already draw one definite conclusion.}

If we assume that \changed{Solution T is the correct solution}, many\footnote{\changed{in the sense of: every scheme that was tested}} popular numerical 
schemes fail to converge to physical solutions.
Although it cannot be ruled out that 
they ultimately converge to Solution T as the numerical grid becomes infinitely fine, they approach Solution N
for computationally accessible grids --- which is all that matters for practical purposes.

On the other hand, if we assume that \changed{Solution N} is the correct solution, 
there is a trivial \emph{theoretical} example of misconvergence: consider the (semidiscrete) Godunov scheme on grids
whose edges are exactly aligned with the discontinuities of \changed{Solution T} (see Figure \ref{fig:example}): in exact arithmetic it would have \changed{Solution T} as steady state on \emph{every} grid.

\changed{In either case --- even if Solution N is correct, which would be less catastrophic for numerical analysis ---} we have to conclude that discrete 
entropy inequalities are not sufficient to avoid convergence to unphysical solutions on feasible grids. Hence they lose a bit of their value as design
principles for numerical schemes, although they are still useful as easy-to-check necessary conditions that are sufficient for scalar conservation laws 
and (probably) 1D systems (as supported by the recent work on small total variation solutions described in Section \ref{section:related-work}).

\changed{Although many reports of deficiencies of various numerical schemes have been published, the clear case of failure observed here has no precedent.}

\if01
The Godunov setup above appears unstable, especially if any amount of artificial viscosity is added, and may be rare in practice. 
However, the well-known ``carbuncle phenomenon'' (see \cite{quirk}), which is rather persistent, is another example of a phenomenon that 
affects especially numerical methods based on Riemann solvers and vanishes if a modicum of artificial viscosity is added to the numerical scheme. 
It is possible that the carbuncle phenomenon and the numerical phenomena discussed in this paper are related; in fact Figure \ref{fig:example} has 
a certain similarity to the numerically observed ``carbuncles''. \cite{robinet-gressier-casalis-moschetta} discovered a new unstable mode in the 
linearized Euler equations for a plane shock wave and conjecture that carbuncles are an inherent instability of the Euler equations 
(TODO: need to be very precise in this quotation), i.e.\ in some
way reflect a particular exact solution. 
Based on the findings in this paper, we propose to modify the conjecture in \cite{robinet-gressier-casalis-moschetta} as follows: 
carbuncles reflect an inherent instability or even nonuniqueness in the Euler equations with EEF condition that does not appear if 
the VV condition is used, in correspondence to the disappearance of carbuncle phenomena if artificial viscosity is added to numerical schemes.
\fi

\section{Theoretical interpretation}

\changed{It remains to discuss which of Solution T and Solution N is 
the physical one and, if Solution T is correct, what causes numerical schemes to produce Solution N.
There are three possible explanations (which are not mutually exclusive):
\begin{enumerate}
\item either Solution N is an example of failure of numerical methods, or
\item EEF solutions are not stable (in the sense of continuous dependence on initial data), or 
\item EEF solutions are not unique.
\end{enumerate}}

\subsection{Breakdown of numerical methods}

\changed{It has already been shown in Section \ref{section:numfail} that the Godunov scheme is flawed, in the sense that it can fail to converge to the physical solution on feasible grids. Hence it is natural to suspect that Solution N is a numerical artifact that does not correspond to a seccond EEF solution (or any other type of solution of the Euler equations).}

\changed{Solution T is steady and self-similar.}
``Steadyness'' is a
non-generic property that is usually not inherited by finite-accuracy numerical solutions
(for example for a Riemann problem that is solved exactly by a single shock, most numerical schemes
produce small additional waves and a slightly different shock).
However, one would expect numerical approximations to be at least \emph{almost} steady, unlike Solution N.

\changed{In a single space dimension, the conservation property of numerical schemes often guarantees accurate shock locations, even if the overall accuracy of the scheme is poor.} 
\changed{On the other hand,}
in two or more dimensions numerical imprecision can significantly change the shape and location of shocks. This may be the cause
of Solution N.
For example, the upwards deflection of the incoming flow by the 
high-pressure area in front of the stagnation region could be weaker in numerical calculations than in Solution T;
the additional pressure would cause the stagnation region to collapse. 

However, in this case the numerical results would depend strongly 
on the choice of \changed{numerical method, mesh width and other parameters.} This is not observed; rather, all choices produce essentially
the same results.

\subsection{Instability}

\changed{A second explanation is instability (in the sense of lack of continuous dependence on the initial data).
It is possible that Solution N is an approximation to an unsteady EEF solution T' that results
from a slight perturbation of solution T at initial time}
(such perturbations are inevitable in most numerical computations 
due to inexact arithmetic, discretization error, artificial viscosity etc.) 
Since \changed{Solution N} is produced (up to minor differences) 
for any ``perturbation'' (i.e.\ for any choice of mesh, numerical method and parameters), 
it would indicate that \changed{Solution T} constitutes a set of initial values for which the Euler equations are not stable.

\changed{On inspection in similarity coordinates, it appears that the numerical solutions are bounded,
converge quickly to Solution N and remain steady 
(many orders of magnitude of time have been observed), so Solution T' would have to be at least approximately self-similar. 
If we assume it is \emph{asymptotically} self-similar for large time, then the following theorem implies that the asymptote
is an EEF solution:
\begin{theorem}
    \label{th:asymptotic}%
    Let $u\in L^\infty(\RC)$ be an EEF solution of \eqref{eq:weak}.
    Assume that $u$ is \defm{asymptotically self-similar} (see Definition \ref{def:asymptotic}),
    then its asymptotic limit $w$ (a self-similar function) is an EEF solution as well.
\end{theorem}
(The proof of Theorem \ref{th:asymptotic} and an analogous result for steady solutions are presented in the Appendix.)
The asymptote would assume the same initial data as Solution T, but would have to be different from it (by closeness to Solution N).
This would already
imply the third explanation (nonuniqueness of EEF solution).} 

\changed{To avoid that, it is necessary to assume either that Solution T' is approximately, but \emph{not asymptotically} self-similar 
--- for example it might oscillate periodically around some self-similar function without approaching it ---, or that
the self-similar asymptote has data at infinity that does not match the initial data (Solution T). Either of these cases
would be revealed by a sufficiently fine numerical grid: the numerical computation on that grid would refuse to
converge to a steady state for the given boundary data.
But the mesh used to compute Solution N is already rather fine, as can be seen from the curved shocks in Figure \ref{fig:second-vx}; 
there is no obvious reason why an even better grid is required.}

\subsection{Nonuniqueness of EEF solutions}

\changed{The third explanation is nonuniqueness: Solution N corresponds to an EEF solution that assumes the same initial data
as Solution T.}

\if\doprivate
%The initial data is reminiscent of expansion shocks (the classical 1D example for unphysical weak solutions and nonuniqueness of weak solutions): 
%in the origin, the gas in the upper and lower post-shock cones moves apart.
%Expansion shocks cannot appear in solutions to (uniformly) viscous perturbations of systems of conservation laws because the smoothing
%effect of viscosity smears the shock
%and the subsequent convection spreads the flow into a rarefaction fan; this happens for arbitrarily small (but positive) viscosity coefficient.
%It is likely that similarly, in Figure \ref{fig:example}, slight smoothing from small viscous terms breaks up the shock pattern in the origin
%and causes the evolution into the self-similar solution in Figure \ref{fig:second-vx}. We believe that Figure \ref{fig:example} is not a VV solution.
\fi

In this context, the following peculiarity is important: according to Proposition
\ref{th:piecewise}, 
the EEF condition \eqref{eq:eef-entropy} is ``insensitive'' to sets with $(d-1)$-dimensional
Hausdorff measure $0$ (such as a single point, for $d=2$);
e.g.\ if \eqref{eq:eef-entropy} is satisfied for $\phi\in C_c^\infty(\overline{\RC}-\{0\})$, it is satisfied for
\emph{all} $\phi\in C_c^\infty(\overline{\RC})$. 
\changed{In verifying that Solution T is an EEF solution, we may ignore the rather singular wedge tip 
in the origin of Figure \ref{fig:example}. It is counterintuitive that conditions for the physical correctness of
solutions may ignore such singularities.}

It seems unlikely that the solution in Figure \ref{fig:example} is stable under small perturbations at the origin 
(such as perturbations from viscous terms in the VV limit).
Note that a planar shock, with inflow state on one side and stagnation area state on the other side, would \emph{not} be steady but move into
the stagnation zone quickly --- it seems unlikely that the example data, which has less mass and energy and more $x$-momentum in the $\{x>0\}$
halfplane, would yield a steady pattern in the origin (however, the ``maximum principle'' implicit in this argument is merely heuristic and
may be wrong in some instances).
%(TODO: why not compute the exact solution of the 1D Riemann problem and compare it to Figure \ref{fig:second} ? )

\if\doprivate
Figure \ref{fig:example} is closely related to the well-known Prandtl-Meyer problem of flow along solid walls with corners.
The physically observed pattern \changed{(which is at the same time} an EEF solution of the Euler
equations) is on the left of Figure \ref{fig:bend} (see Section 111 in \cite{courant-friedrichs} for discussion). 
On the right is another steady EEF solution; as in Figure \ref{fig:example}, there is flow along a contact discontinuity 
with a stagnation area on the other side. This solution is not observed in physical\footnote{\changed{It is not observed in sophisticated numerical computations either. However, some classical schemes --- particularly the Godunov scheme --- are prone to produce the solution on the right. This adds another example for numerical failure (similar to Section \ref{section:numfail}
and for insufficiency of the EEF condition. However, this example is not as strong as the one in Figures \ref{fig:example} and F\ref{fig:second-vx} because it has a solid boundary.}} experiments.
This analogy is further evidence that Figure \ref{fig:second-vx} is indeed an approximation to the physically correct solution for the initial
data in Figure \ref{fig:example}.

\begin{figure}
\begin{center}
\input{bend.pdftex_t}
\end{center}
\caption{Supersonic inflow from the left reaches a corner. Left: physical solution. Right: solution with a stagnation area behind the corner;
this solution is not observed in experiments.}
\label{fig:bend}
\end{figure}
\fi

But it is the following observation that provides the strongest argument for nonuniqueness: 
the Lax-Wendroff theorem (see \cite{lax-wendroff}; see also \cite{godlewski-raviart}, \cite{kroener-rokyta-wierse} and
most generally \cite{elling-lax-wendroff} for Lax-Wendroff-type theorems for irregular grids) states that if a numerical scheme is consistent
and satisfies a discrete entropy inequality (see \cite{harten-hyman-lax,majda-osher-discrete-entropy,osher-chakravarthy,tadmor-entropy-mathcomp,tadmor-entropy-stable,osher-tadmor}), then
the limit of a boundedly almost everywhere converging sequence of numerical solutions is an EEF solution.
\changed{(Note that we prescribe the full numerical flux on the boundaries, so the boundary conditions are analogous to an initial condition and can be treated with a straightforward modification of the Lax-Wendroff theorem for initial-value problems.)} 
The Godunov scheme, used to compute Figure \ref{fig:second-vx}, is
the standard example for a consistent scheme that satisfies all discrete entropy inequalities.
Our numerical solutions do, on inspection, appear to converge quickly; this would imply that Solution N corresponds to an EEF solution.

\section{Related work}

\label{section:related-work}

For multidimensional scalar ($m=1$) conservation laws with arbitrary $f$, 
\cite{kruzkov} (generalizing earlier work) shows that a global EEF solution exists, is unique, satisfies the VV condition as well,
and is stable under $L^1$ perturbations of the initial data.

\cite{glimm} provides a famous existence proof for strictly 
hyperbolic \emph{systems} with genuinely nonlinear fields 
and initial data with small total variation; the interaction functionals constructed in this paper
are a crucial ingredient for all subsequent work. \cite{liu-admissibility-memoir} 
extends the result to systems with linearly and some nonlinearly degenerate 
fields. \cite{bressan-crasta-piccoli} constructed the Standard Riemann Semigroup (SRS), an $L^1$-stable semigroup
of EEF solutions for initial data with small total variation, for strictly hyperbolic systems 
with genuinely nonlinear or linearly degenerate fields (see also \cite{liu-yang}).
\cite{bressan-lefloch} showed that EEF solutions to 1D systems are unique and coincide with the SRS solutions, 
under certain smoothness assumptions including small total variation (see also \cite{bressan-goatin}). 
\cite{bianchini-bressan} prove that for small $TV$ initial data and strictly 
hyperbolic (but otherwise arbitrary) systems VV solutions exist and are stable under $L^1$ perturbations of the initial data,
so for the class of solutions that are subject both to \cite{bressan-lefloch} and to \cite{bianchini-bressan}, EEF and VV solutions 
are equivalent.

On the other hand, an EEF pair $(\eta,\Vec\psi)$ has to satisfy the condition \eqref{eq:eef} which is an overdetermined problem for $m\geq 3$, 
so for some systems no EEF pairs exist and the EEF condition is void. However, EEF pairs do exist
for most physically relevant systems, even those with $m\geq 3$. More seriously, for certain $2\times 2$ systems (with nonlinear degenerate fields) 
\cite{conlon-liu} construct a single weak shock that is an EEF solution but does not satisfy the Liu entropy condition 
(see \cite{liu-condition-uniqueness-twobytwo,liu-condition-uniqueness-general}). By \cite{bianchini-bressan}, there must be a VV solution (for
the same initial data) that satisfies the Liu entropy condition as well --- so it cannot be the aforementioned weak shock. Therefore
the example in \cite{conlon-liu} also constitutes an example of a nonunique EEF solution, albeit for an ``artificial'' system with nonlinear
degeneracy.

\cite{hopf-eef} proposes the EEF condition for scalar conservation laws ($m=1$), proves that it is implied by the VV condition under some circumstances
and notes that there is a large set of convex entropies. Apparently independently, 
\cite{kruzkov} obtained analogous results for systems. \cite{lax-zarantonello} contains the first use of the term ``entropy condition'' 
for the EEF condition. 
Various forms of the EEF condition had been known and in use for special systems such as the Euler equations for a long time 
(e.g.\ by the name of \defm{Clausius-Duhem inequality}), 
especially as shock relations; however, the above references seem to be the first to define the general notion of strictly convex EEF pairs, 
to propose the EEF condition as a mathematical tool for arbitrary systems of conservation laws and to formulate it in the weak
form \eqref{eq:eef-entropy} rather than the special case \eqref{RH-EEF}.

\cite{li-zhang-yang} provide an analytical and numerical discussion of 2D Riemann problems for various systems including the 
Euler equations. However, they focus on data constant in each of the four quadrants, so Solution T is not covered.

\if01
(TODO: what about viscous traveling wave solutions for large TV data riemann problems for Euler equations?)
(TODO: what about $2\times 2$ conservation laws with large data?)
(TODO: look at diperna paper cited in goatin/lefloch article.)

(TODO: mention that \cite{courant-friedrichs} fig 30 p. 296 is our example if the wedge is replaced by stagnation air;
see p. 345 Fig 69. 
Quote:
``All these and other mathematically possible flow patterns with a singular center Z are at our disposal for interpreting
experimental evidence. Which, if any, of these possibilities occurs under given circumstances is a question that cannot
possibly be decided within the framework of a theory with such a high degree of indeterminacy. Here we have a typical instance
of a theory incomplete and oversimplified in its basic assumptions; only by going more deeply into the physical basis of our theory,
i.e.\ by accounting for heat conduction and viscosity, can we hope to clarify completely the phenomena at a three-shock singularity.
It may well be that the boundary layer which develops along the constant discontinuity line modifies the flow pattern sufficiently
to account for the observed deviation; [...quote Liepmann paper]''
\fi

\section{Conjectures and final remarks}

\if\doprivate
The total variation of the initial data cannot be made arbitrarily small because the oblique shock relations can be solved only
for inflow velocity above some supersonic limit (which depends on $\alpha$). 
%%It is possible to make all density jumps 
%%and the velocity jump across the shock small by choosing weak shocks and small $\alpha$ (but the contacts retain a strong velocity jump), or to 
%%make the velocity jump across the contacts weak at the expense of using strong shocks. 
It would be interesting to find modified examples with
arbitrarily small total variation.
\fi

\noindent The results demonstrate that

\emph{the Godunov method does not always converge to the physical solution on feasible grids.}

%%Assuming that EEF solutions are unique (which the author does not believe any more), they would demonstrate that
%%\emph{many popular numerical methods fail to converge to the correct solution (on feasible grids).}

\noindent Moreover, they suggest the following conjecture:

\emph{EEF solutions to the multidimensional Euler equations are not always unique.}

\noindent If this conjecture is true, it would have far-reaching consequences. 
The EEF condition would not be sufficient as a selection principle for physical/unique solutions, except in special cases like the 
ones described in Section \ref{section:related-work}. It would be necessary to find ways to use the cumbersome VV condition or to discover new entropy conditions. 

Although the numerical results support the conjecture unambiguously, the question is so important that a rigorous proof is highly desirable. 
However, since the initial data has large vorticity at the contact discontinuity, it seems difficult to construct (or to prove results about) 
exact solutions. One possible line of attack is to derive novel entropy conditions from the VV condition and to check whether they are violated 
by the steady solution in Figure \ref{fig:example}.

In any case, this paper motivates the investigation of multidimensional Riemann problems for systems; these appear to be very difficult 
and exhibit a large variety of phenomena (see \cite{lax-xdliu,li-zhang-yang}). 
This goal requires techniques for proving existence of smooth steady or self-similar solutions to boundary-value problems for systems 
of nonlinear hyperbolic conservation laws; while there are classical methods for smooth solutions in hyperbolic regions, 
work on tools for the elliptic and mixed case has begun only recently (see e.g.\ \cite{elling-liu-ellipticity-journal}).
% what about serre?

\if\doprivate
Finally, one could wonder whether Figure \ref{fig:example} is an example of a ``generic'' phenomenon or whether EEF solutions are unique for 
all initial data outside a ``small'' complement. This question is being investigated; we suspect that nonuniqueness
is indeed generic.
\fi

\section*{Appendix: asymptotically steady and self-similar weak solutions}
\par

Remark: in Theorem \ref{th:asymptotic} and in the following statements, 
\[ f(t,\cdot) \rightarrow g \qquad\text{in $L^1_{\loc}(\Omega)$} \]
as $t\downarrow 0$ resp.\ $t\uparrow\infty$ 
is to be understood as: for all $\epsilon>0$ and $K\Subset\Omega$ there is a $T=T(\epsilon)>0$ so that
for almost all $0<t\leq T$ resp.\ $t\geq T$,
\[ \| f(t,\cdot)-g \|_{L^1(K)} \leq\epsilon. \]

\begin{lemma}
    \label{lemma:weak}%
    Let $u\in L^\infty(\RC)$, $u_0\in L^\infty(\Rd)$. If
    \begin{enumerate}
    \item
        $u(t,\cdot)\rightarrow u_0$ in $L^1_{\loc}(\R^d)$ as $t\downarrow 0$, and
    \item
        $u$ satisfies \eqref{eq:weak} for all $\phi\in C^\infty_c(\underline{(}0,\infty)\times\R^d)$,
    \end{enumerate}
    then $u$ satisfies \eqref{eq:weak} for all $\phi\in C^\infty_c(\underline{[}0,\infty)\times\R^d)$.
\end{lemma}

\begin{proof}
    Let $\theta\in C^\infty[0,\infty)$ so that $\theta=1$ on $[0,1]$, $\theta=0$ on $[2,\infty)$. For any $T>0$, define
    $\theta_T(t):=\theta(T^{-1}t)$. Note that $|\theta_T|=O(1)$, $|\theta_T'|=O(T^{-1})$ (as $T\downarrow 0$). For given $\phi\in C_c^\infty([0,\infty)\times\Rd)$,
    split $\phi_1(t,x):=\theta_T(t)\phi(t,x)$ and $\phi_2=\phi-\phi_1$.
    \begin{alignat*}{1}
        &\int_0^\infty\int_{\Rd}u\phi_t+\Vec f(u)\cdot\nabla\phi~dx~dt \\
        &= \int_0^{2T}\int_{\Rd}u\phi_{1t}+\Vec f(u)\cdot\nabla\phi_1~dx~dt + \int_T^\infty\int_{\Rd}u\phi_{2t}+\Vec f(u)\cdot\nabla\phi_2~dx~dt.
    \end{alignat*}
    Since $\phi_2\in C_c([T,\infty)\times\Rd)$, the second summand vanishes by assumption. The first summand equals
    \begin{alignat*}{1}
        &= O\left(T\sup_{0<t\leq 2T}\|u(t,\cdot)-u_0\|_1\cdot(T^{-1}+1)\right) 
        + \int_0^{2T}\int_{\Rd}u_0\phi_{1t}~dx~dt\\
        &+O\left(T(\|u_0\|_1+\sup_{0<t\leq 2T}\|u(t,\cdot)-u_0\|_1\right) \\
        &= O\left(\sup_{0<t\leq 2T}\|u(t,\cdot)-u_0\|_1\right) - \int_{\Rd}u_0(x)\phi(0,x)~dt \\
        &+ O\left(T(\|u_0\|_1+\sup_{0<t\leq 2T}\|u(t,\cdot)-u_0\|_1\right).
    \end{alignat*}
    On taking $T\downarrow 0$, all $O$ terms vanish; hence $u$ satisfies \eqref{eq:weak}.
\end{proof}

Remark: the converse of Lemma \ref{lemma:weak} (which is not needed) is not immediate because $u(t,\cdot)-u_0$ may be large for some $t$ as long
as the set of such $t$ has small measure near $0$.

\begin{definition}
    \label{def:asymptotic}%
    \begin{enumerate}
    \item
        A function $u\in L^1_{\loc}(\RC;\R^m)$ is called \defm{asymptotically self-similar} if there is a function $w:\R^d\rightarrow\R^m$ 
        so that 
        \[ u(t,t^{-1}\cdot) \rightarrow w\qquad\text{in $L^1_{\loc}(R^d)$}. \]
    \item
        $u$ is called \defm{self-similar} if, for some $w$, $u(t,t^{-1}\cdot)=w$ for almost all $t>0$.
    \end{enumerate}
\end{definition}

\begin{proof}(of Theorem \ref{th:asymptotic})
    By Lemma \ref{lemma:weak}, to show that $w$ is a weak solution it is sufficient to check that
    \begin{alignat}{1}
        \int_0^\infty\int_{\R^d}w\left(\frac{x}{t}\right)\phi_t(t,x)+\Vec f\left(w\left(\frac{x}{t}\right)\right)\cdot\nabla\phi(t,x)~dx~dt 
        &= 0 \label{eq:asymptotic-1}
    \end{alignat}
    for all $\phi\in C^\infty_c(]0,\infty))$. The essential idea is to scale coordinates to shift the support of $\phi$ into a large-$t$ region and
    to use asymptotic convergence.

    Let $0<t_1<t_2$ be such that $\supp\phi\subset[t_1,t_2]\times\R^d$. Let $\epsilon>0$ be arbitrary, set $T=T(\epsilon)$ as in Definition \ref{def:asymptotic}.
    The change of coordinates $t=\frac{t_1}{T}\tau$, $x=\frac{t_1}{T}\xi$ changes the left-hand side of \eqref{eq:asymptotic-1} into
    \begin{alignat}{1}
        & \left(\frac{t_1}{T}\right)^{d+1}\int_T^{\frac{t_2T}{t_1}}\int_{\Rd}w\left(\frac{\xi}{\tau}\right)\phi_t\left(\frac{t_1}{T}\tau,\frac{t_1}{T}\xi\right)
        + \Vec f\left(w\left(\frac{\xi}{\tau}\right)\right)\cdot\nabla\phi(\tau,\xi)~d\xi~d\tau \notag \\
        & = \left(\frac{t_1}{T}\right)^{d+1}\int_T^{\frac{t_2T}{t_1}}\int_{\Rd}u(\tau,\xi)\phi_t\left(\frac{t_1}{T}\tau,\frac{t_1}{T}\xi\right)
        + \Vec f(u(\tau,\xi))\cdot\nabla_x\phi\left(\frac{t_1}{T}\tau,\frac{t_1}{T}\xi\right)~d\xi~d\tau \notag \\
        & + O\left(\left(\frac{t_1}{T}\right)^{d+1}\cdot T\cdot\epsilon T^d\right) \label{eq:asymptotic-2}
    \end{alignat}
    where $O$ is with respect to $\epsilon\rightarrow\infty$.
    Note that the support of the scaled $\phi$ is in $[T,\infty)\times\R^d$. 
    Also, the assumption that $u$ is bounded is essential here. The second summand on the right-hand side equals
    \begin{alignat*}{1}
        & \int_T^{\frac{t_2T}{t_1}}\int_{\Rd}u(\tau,\xi)\frac{T}{t_1}\phi\left(\frac{t_1}{T}\tau,\frac{t_1}{T}\xi\right)_\tau + 
        \Vec f(u(\tau,\xi))\cdot\nabla_\xi(\phi(\tau,\xi))~d\xi~d\tau.
    \end{alignat*}
    Since $u$ is assumed to be a weak solution, this term vanishes. Taking $\epsilon\downarrow 0$ in \eqref{eq:asymptotic-2} yields
    \eqref{eq:asymptotic-1}.
    
    For the proof of the EEF part, replace $u,w$ by $\eta(u),\eta(w)$ and $f(u),f(w)$ by $\psi(u),\psi(w)$ above.
\end{proof}

The same results as for self-similar weak solutions can be obtained for steady solutions:

\begin{definition}
    \label{def:asymptotically-steady}%
    \begin{enumerate}
    \item
        $u\in L^1_{\loc}(\RC;\R^m)$ is called \defm{steady} if, for some $w:\R^d\rightarrow\R^m$, $u(t,\cdot)=w$ for almost all $t>0$.
    \item
        $u\in L^1_{\loc}(\RC;\R^m)$ is called \defm{asymptotically steady} if there is a $w:\R^d\rightarrow\R^m$,
        so that 
        \[ u(t,\cdot) \rightarrow w\qquad\text{in $L^1_{\loc}(\R^d)$} \]
        for almost all $t\geq T$.
    \end{enumerate}
\end{definition}

\begin{theorem}
    If $u\in L^\infty(\RC;\R^m)$ is an asymptotically steady and bounded weak solution, then $w$ (as in Definition \ref{def:asymptotically-steady}) 
    is a weak solution as well.
    If $u$ is an EEF solution, so is $w$.
\end{theorem}
\begin{proof}
    Let $\phi\in C^\infty_c((0,\infty)\times\Rd)$ be arbitrary. Let $\supp\phi\subset[0,\tau]$.
    For any $\epsilon>0$,
    \begin{alignat*}{1}
        &  \int_0^\infty\int_{\Rd}w(x)\phi_t(t,x)+\Vec f(w(x))\cdot\nabla\phi(t,x)~dx~dt \\
        &= \int_0^\infty\int_{\Rd}u(T+t,x)\phi_t(t,x)+\Vec f(u(T+t,x))\cdot\nabla\phi(t,x)~dx~dt + O(\tau\epsilon\|D\phi\|_\infty) \\
        &= O(\tau\epsilon\|D\phi\|_\infty) 
    \end{alignat*}
    because we can extend $\phi(\cdot-T,\cdot)\in C_c((T,\infty)\times\Rd)$ smoothly by $0$ to a map $\tilde\phi\in C_c((0,\infty)\times\Rd)$
    and use that $u$ is a weak solution. Lemma \ref{lemma:weak} shows that $w$ is a weak solution.

    For the proof of the EEF part, replace $u,w$ by $\eta(u),\eta(w)$ and $f(u),f(w)$ by $\psi(u),\psi(w)$ above.
\end{proof}

\section*{Acknowledgements}

The author would like to thank Tai-Ping Liu for his support and comments and 
Ron Fedkiw, Doron Levy, Wolfgang Dahmen and Ralf Mass\-jung for sharing their insight about numerical schemes.

\if01
\section*{More speculation}

Colella and Woodward say that Godunov produces an unphysical shock at the forward-facing step. But then we know that Godunov should not produce
unphysical things. The truth is that it does, it's just that physical is not the same as entropy-admissible: their example is precisely
the Prandtl-Meyer corner where one expects a rarefaction wave, but a contact (or shock here, since waves impinge on the contact) solution is also possible.

\fi

\bibliographystyle{amsalpha}
\bibliography{../../mypapers/thesis/elling}
%\vskip2cm
%\noindent Preprint version pre-8 (\number\month/\number\day/\number\year)

\end{document}